\documentclass[preprint,12pt]{elsarticle}

\usepackage{amssymb}

\usepackage{moreverb}
\usepackage{amsmath,bm}
\usepackage{subfig}
\usepackage{graphicx}
\usepackage{comment,physics}
\usepackage{stackengine}
\usepackage{lineno}
\usepackage{xcolor}
\usepackage{array,multirow}
\usepackage{algorithm,algorithmic,ulem}

\newtheorem{theorem}{Theorem}

\newlength\myindent
\newcommand\bindent{%
	\begingroup
	\setlength{\itemindent}{\myindent}
	\addtolength{\algorithmicindent}{\myindent}}
\newcommand\eindent{\endgroup}

\journal{Journal of Computational Physics}

\begin{document}

\begin{frontmatter}



\title{Laplace-Beltrami based Multi-Resolution Shape Reconstruction on Subdivision Surfaces}


\author[1]{A. M. A. Alsnayyan\corref{cor1}}
\ead{alsnayy1@msu.edu}
\author[1]{B. Shanker}

\address[1]{Electrical and Computer Engineering, Michigan State University, Michigan, USA}

\cortext[cor1]{Corresponding author}

\begin{abstract}
The eigenfunctions of the Laplace-Beltrami operator have widespread applications in a number of disciplines of engineering, computer vision/graphics, machine learning, etc. These eigenfunctions or manifold harmonics, provide the means to smoothly interpolate data on a manifold. They are highly effective, specifically as it relates to geometry representation and editing; manifold harmonics form a  natural basis for multi-resolution representation (and editing) of complex surfaces and functioned defined therein. In this paper, we seek to develop the framework to exploit the benefits of manifold harmonics for shape reconstruction. To this end, we develop a highly compressible, multi-resolution shape reconstruction scheme using manifold harmonics. The method relies on subdivision basis sets to construct both boundary element isogeometric methods for analysis and surface finite elements to construct manifold harmonics. We pair this technique with volumetric source reconstruction method to determine an initial starting point.  Examples  presented highlight efficacy of the approach in the presence of noisy data, including significant reduction in the number of degrees of freedom for complex objects,  accuracy of reconstruction, and multi-resolution capabilities.

\end{abstract}




\end{frontmatter}


\section{Introduction}

Research into inverse scattering dates back decades and has found applications in a number of wide ranging fields of studies, including areas such as 
medical diagnostics, detection of buried objects, tomography, and non-destructive evaluation \cite{EM_buried,diff_tomo,optical_tomo_optim,ultrasound_survey,subsurface,Langenberg1993}; in these problems, the goal is to retrieve the distribution of constitutive properties in a domain and/or geometry given a set of measured scattered field data. Along these lines of an inverse scattering problem, there are two subclasses of problems that can be considered: can one (a)  modify shapes such that one obtains desired scattered fields, or (b) reconstruct the shape of a object given scattered field data (and boundary conditions on the surface). As is evident, both problems are closely related. The former is commonly referred to as shape optimization; it has a number of applications ranging from acoustics \cite{acoustic_horn,cirak_fem_shape,CAD_optim,structural_acoustic,kress_phaseless} to electromagnetics \cite{EM_SHAPE_OPTIMZATION,RCS_optimization} to  medical imagining \cite{optical_tomo_optim,ultrasound_survey} to design of horns \cite{acoustic_horn}, and a number of other applications \cite{Jan_multiresolution,Aerodynamic_Optimization}. The computational techniques used therein have been integrated with gradient based optimization methods \cite{cirak_subd_statics,RCS_optimization,kress_phaseless} as well global optimization methods \cite{mich_GA_EM,PSO}. 

This paper is devoted to the latter problem - the recovery of scattering obstacles from phaseless far-field data. Interest in this class of problems is widespread as it finds application in a number of different disciplines \cite{kress_colton_book}. The main challenges that arise in this problem are (a) the ill-posed nature of the problem, (b) the number of degrees of freedom in the optimization problem, and (c) the optimal minimization method used for optimization. In this paper, the aim is to address challenge (b) in a shape reconstruction optimization routine that is constrained by the forward scattering problem with the goal of minimizing the geometrical parameters that define the boundary of the scatterer, or in other words minimizing the search space. In particular, we look at constructing and modifying a compressed parameterization of some starting surface until we reach some minimization of a cost functional that measures the discrepancy between the prescribed far field data and the far field pattern corresponding to the current approximation of the paramterized scatterer. As is apparent, the geometry parameterisation and its interplay with the analysis plays a crucial role in shape optimisation/reconstruction. When a geometric representation, distinct from the computer-aided design (CAD) model is used to represent the shape, additional errors are incurred as well as a large parameterisation space is needed to represent high fidelity shapes, drastically increasing the difficulty of the optimization scheme. In addition, it leads to non-physical oscillations in the optimised geometry [17,18] as well as severely distorted mesh requiring auxiliary mesh smoothing, further complicating the problem. To remedy these difficulties, the shape optimisation/reconstruction can be constructed in a high-order isogeometric framework. 

Isogeometric methods use the same underlying basis sets to represent both the geometry and physics on the geometry. This class of techniques was pioneered by Hughes \cite{hughes2005isogeometric}, and has since been adapted for a number of different types of problems in structural mechanics \cite{cirak_fem_shape,CAD_optim}, electromagnetics \cite{jiesimply,Jiescalar,JieDaultShankerChapterSubd,Fu2017GeneralizedDS} and acoustics \cite{abdel_acoustics,acoustics_IGA_cirak}. The research in using isogeometric analysis has treaded along two paths; use of non-unifrom B-splines and subdivision surfaces; our focus is on subdivision. Subdivision is a powerful geometric modelling technique for generating smooth surfaces on arbitrary connectivity meshes which are the generalisation of splines to arbitrary connectivity meshes. Specifically, we use the Loop scheme based on triangular meshes and quartic box-splines \cite{loop1987smooth}. Subdivision surfaces are perhaps the ideal candidate for IGA shape reconstruction as one can exploit other facets of subdivision--hierarchical refinement, higher order continuity, and arbitrary topology, i.e., no other restricting assumptions on the geometry need be considered (star-shaped domains, for example), see Refs.~\cite{iga_shape,cirak_fem_shape,efficient_shape} and references therein for examples on IGA shape optimization. The general approach in these papers is to exploit the hierarchy of the control mesh underlying a subdivision surface to do multiresolution editing \cite{multires,mesh-edit,cirak_subd_statics}. The coarse control mesh vertex positions are modified to perform large-scale editing and the fine control mesh vertex positions are modified to add localized changes. This allows for a relatively compressed representation of geometry, allowing for a smaller search space and thus computationally feasible shape reconstruction problem. This approach works well. But, one of the questions that we ask in the paper is \emph{whether} we can construct a better compression scheme. An inspiration for the answer to this question can be found in computer graphics \cite{LBO_DNA,LBO_MH,Italian_LBO_book,LBO_UND_GEOM}. The nub of these ideas is to develop a mechanism that enables compression and morphing of manifolds in an efficient manner.  

These techniques rely on the manifold harmonic basis (MHB) set; a basis set constructed from the eigenfunctions of the Laplace-Beltrami operator (LBO). As is well known, LBO is a  self-adjoint linear surface operator that captures all the intrinsic properties of the shape and is invariant to extrinsic shape transformations such as isometric deformations. Its eigenfunctions or the MHBs can be understood as a generalization of the Fourier spectrum for functions defined on a general surface manifold. They provide a unique, compact, elegant, and multi-resolution basis for spectral \emph{shape} processing that is independent of the actual shape representation and parameterization. Several successful applications have been proposed that take advantage of these desirable properties, such as spectral geometry filtering, compression, and surface deformation \cite{LBO_MH,Poly_mesh_proc}. 

In this paper, we will leverage our earlier work on subdivision-based isogeometric methods for acoustics to develop a framework that carries over the benefits of subdivision meanwhile maintaining the favorable properties of MHs for use in both analysis and morphing the manifold for shape reconstruction. The specific contributions of this paper are as follow:
\begin{enumerate}
	\item develop a MHB based compression scheme for representation of and analysis on manifolds, 
	\item develop an inverse source based initial estimate for shape, 
	\item develop a multi-resolution framework for shape reconstruction, 
	\item and demonstrate the viability of this technique on a number of challenging targets. 
\end{enumerate}
The methods developed in this paper are agnostic to the specific technique used in the optimization procedure used while reconstruction. In other words, the novelty/contributions of this paper do not lie in the optimization scheme chosen, but the multiresolution framework used to systematically update and refine the geometry and physics throughout the optimization process. And as such, we have used readily available libraries for the optimization procedure (Method of Moving Asymptotes (MMA) \cite{NLopt}) and use the straight-forward finite difference method to implement the optimization. If one were to use a gradient based approach, there exist more efficient adjoint based methods that can be implemented \cite{KRESS199249,Kirsch_1993,Frechet_diff,optical_tomo_optim,LeLouer2018}.

The remainder of this paper is organized as follows: In Section 2, the forward scattering constraint and the shape reconstruction problem are posed. Sections 3 and 4 elucidates the details the Loop subdivision boundary element formulation of the forward problem. In Section 5, we formulate the MHBs, demonstrate its spectral properties, and develop a backprojection scheme to initialize our shape reconstruction.  In Section 6, we present a number of results on structurally challenging objects to demonstrate the salient features of our reconstruction algorithm. Finally, in Section 7, we summarize our contribution in this body of work.

\section{The forward scattering and shape reconstruction problems}
\begin{figure}[!tbh]
	\centering
	\includegraphics[width=0.2\linewidth]{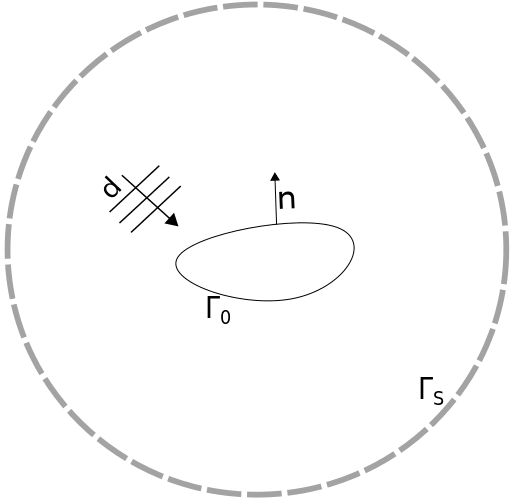}
	\caption{The shape reconstruction problem.}
	\label{inverse_setup}
\end{figure}

Consider the problem depicted in Figure \ref{inverse_setup}: Here, there exists a soft scatterer embedded in $\Omega \in \mathbb{R}^3$, whose boundary is denoted by $\Gamma_0$ with a uniquely defined outward pointing normal $\hat{\vb{n}}$. Without loss of generality, assume that there exists a spherical surface $\Gamma_s \in \mathbb{R}^3$ that circumscribes $\Gamma_0$, but is sufficiently removed from it.

On $\Gamma_s$, data due to fields scattered by the obstacle, $\Phi_{G}^{s}(\vb{r})_{\left ( \vb{d}_{i},\kappa_{j} \right )}$, is available for a combination of incident waves propagating in direction $ \vb{d}_{i} \in \mathbf{D}$, where $\mathbf{D} = \left \{ \mathbf{d}_1, \mathbf{d}_2, \cdots, \mathbf{d}_n \right \} $, with wavenumbers $\kappa_j \in \vb{K}$, where  $\vb{K} = \left \{ \kappa_1, \kappa_2, \cdots, \kappa_k \right \}$, denoted as $\Phi^i (\vb{r})_{\left ( \vb{d}_{i},\kappa_{j} \right )}$ (in the what follows, the subscripts will be omitted for notational simplicity). For each incident field, the resulting total field $\Phi^t (\vb{r}) = \Phi^{s} (\vb{r}) + \Phi^{i} (\vb{r}) $ satisfies the following boundary value problem
\begin{subequations}
	\begin{align}
	\nabla \Phi^{t}(\mathbf{r}) + k^{2}\Phi^{t}(\mathbf{r}) &= 0 \hspace{1cm} \textbf{r} \in \Omega     \label{eq:helmholtz_equation}, \\
	\Phi^{t}(\mathbf{r}) &= 0 \hspace{1cm} \textbf{r} \in \Gamma_{0},     \label{eq:boundary_condition} \\
	\lim_{r\rightarrow \infty}\sqrt{r}\left (\frac{\partial \Phi^{s}}{\partial n} - i\kappa \Phi^{s} \right ) &= 0  \hspace{1cm} \textbf{r} \in \Omega.         \label{eq:rad_condition}
	\end{align}
	\label{eq:forward_problem}
\end{subequations}
Here, $\Phi^{s}$ denotes the scattered field, \eqref{eq:boundary_condition} defines the sound-soft boundary condition and \eqref{eq:rad_condition} imposes the Sommerfeld radiation condition. 

The objective of this experiment is to determine the shape of the scatter $\Gamma_{0}$ such that the measured scattered fields $\Phi^{s}_{m}(\vb{r})$ at $\Gamma_{s}$ obtained for a collection of directions and wavenumber $\left( \vb{d}_{i},\kappa_{j} \right)$ satisfy $\Phi^{s}_{m} (\vb{r}) = \Phi_{G}^{s}(\vb{r})$. This is tantamount to solving a minimization problem given by
\begin{equation}
\label{eq:objective_function}
\text{find} \mathop{argmin}_{\Gamma \in \mathbb{R}^{3}}J(\Gamma) := \frac{1}{2} \int_{\Gamma_{s}}| \Phi_{G}^{s}(\vb{r}) - \Phi^{s}_{m}(\vb{r})|^{2}d\vb{r}. 
\end{equation}
Note, as is evident from the objective function, we assume that phaseless data is available on the surface $\Gamma_s$. As stated, this is a constrained shape reconstruction problem with $\Gamma$ as the design variable and the forward scattering problem as the constraint on the admissible exterior fields $\Phi^{s}_{m}$. In what follows, we use far-field scattered data, in which case the objective function is defined in the limiting case in which $\Gamma_{s}$ goes to infinity; following standard practice, fields are normalized by distance. 

\begin{figure}[!bh]
	\centering
	\includegraphics[width=0.4\linewidth]{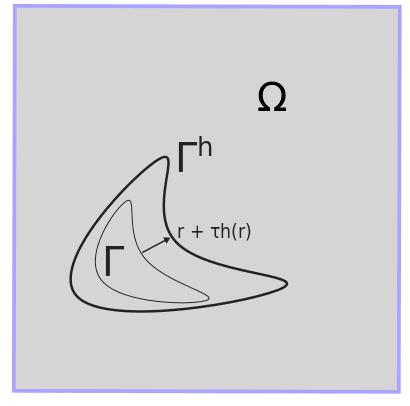}
	\caption{The perturbed shape according to design
		perturbation \textbf{h}(\textbf{r}).}
	\label{fig:perturbed_surface}
\end{figure}

The strategy we use to minimize the cost functional (\ref{eq:objective_function}) is an iterative gradient-based approach which requires the derivative of the cost functional (\ref{eq:objective_function}) with respect to the domain perturbations. The domain perturbation, defined by $\Gamma \rightarrow  \Gamma^{h_{\tau}} $ is described by a vector displacement field $\textbf{h}(\textbf{r})$ such that $\textbf{r} \rightarrow \textbf{r} + \tau \textbf{h}(\textbf{r})$, where
$\textbf{r} \in \Gamma$ and
$\tau$ is a scalar parameter that denotes the amount of shape change, as illustrated in Fig. \ref{fig:perturbed_surface}. The reader is referred to \cite{Jasbir_senstivity,intro_shape}, and references therein for a more detailed definition of shape perturbation as we have only outlined a brief sketch here. To this end, we define the gradient of the cost functional (\ref{eq:objective_function}) as
\begin{equation}
\label{eq:der_obj_func}
J^{\prime}(\Gamma)  =  \lim_{\tau \rightarrow 0} \frac{J(\Gamma^{h_{\tau}}) - J(\Gamma)}{h}.
\end{equation}
Where $J(\Gamma)$ and $J(\Gamma^{h_{\tau}})$ are the cost functional evaluated for the reference and perturbed domain.

\subsection{The objective function}
In our shape reconstruction scheme, the evaluation of the objective function and the gradient of the objective function are some of the key elements needed to minimize the objective function. To this end, the shape of the boundary and the perturbation vector field \textbf{h} is described by a linear combination of manifold harmonics (detailed in later sections). As a result, we have entailed the boundary in terms of shape coefficients $\beta = \left\{\beta_{i}\right\}$, wherein $\Gamma(\beta)$. We can now repose the reconstruction problem as
\begin{equation}
\label{eq:objective_function_MHs}
\text{find } \beta  \text{ such that } J(\Gamma(\beta)) := \mathop{min}\frac{1}{2} \int_{\Gamma_{s}}| \Phi_{G}^{s}(\vb{r}) - \mathcal{L}_{far}[\Lambda,\Gamma\left(\beta \right)]|^{2}d\vb{r}. 
\end{equation}
where $\mathcal{L}_{far}$ is the operator that maps the boundary data $\Lambda = \frac{\partial \Phi}{\partial n}$ onto the farfield for some object $\Gamma$, described by a set of shape coefficients $\beta$. The gradient of the objective function $J$ with respect to the geometrical parameter $\beta_{i}$ can be done using a finite difference method \cite{vogel2002computational}
\begin{equation}
\frac{\partial J(\beta)}{\partial \beta_{i}} = \frac{J(\beta_{1},\cdots,\beta_{i}+\tau_{i},\cdots,\beta_{N_{v}}) - J(\beta_{1},\cdots,\beta_{i},\cdots,\beta_{N_{v}})   }{\tau_{i}}.
\end{equation}
This finite difference scheme is the most straightforward implementation, given that it still enables us to highlight the main objective of this paper: the use of multiresolution MHs in shape optimization. The reader is turned to \cite{optical_tomo_optim,vogel2002computational} for more efficient approaches that can be adopted to this scheme. The evaluation of $\mathcal{L}_{far}$ in effect, the objective function $J$, for some shape is done by discretizing our forward problem \eqref{eq:forward_problem} by the Boundary Element Method, which is detailed in the following section.

\section{Boundary Element Formulation}

In this section, we will provide a general outline on applying  boundary element method to solve (\ref{eq:forward_problem}) and obtain the measured far-field scattered fields $\Phi^{s}_{m}(\textbf{r})$, where $\textbf{r} \in \Gamma_{s}$, in order to compute the objective function and its gradient. In particular, the Burton-Miller formulation \cite{BurtonMiller} provides a unique solution, and is expressed as: 
\begin{equation}
\begin{split}
\mathcal{L}_{BM} \left [ \Lambda,\Gamma\right] (\mathbf{r})& = V^{i} (\mathbf{r})\hspace{4.7cm} \textbf{r} \in \Gamma,\\ 
\mathcal{L}_{BM} \left [ \Lambda,\Gamma\right](\mathbf{r}) & \doteq  (1-\alpha) \mathcal{S}[\Lambda,\Gamma](\mathbf{r})  + \alpha \beta \mathcal{D}'[\Lambda,\Gamma](\mathbf{r}) \hspace{0.41 cm} \textbf{r} \in \Gamma,\\ 
V^{i} (\mathbf{r}) & \doteq (1-\alpha)\Phi^{i}(\mathbf{r}) + \alpha \beta \hat{\textbf{n}}\cdot \nabla \Phi^{i}(\mathbf{r}) \hspace{1.2 cm} \textbf{r} \in \Gamma.
\end{split}
\label{eq:BurtonMiller}
\end{equation}
where $\Lambda(\mathbf{r})$ denotes $\partial_\mathbf{n} \Phi (\mathbf{r})$,  $\hat{\textbf{n}}$ is the outward unit normal vector to $\Gamma$, $0 \leq \alpha \leq 1$ is a coupling factor that makes the solution unique at all frequencies \cite{BurtonMiller}, and $\beta$ is the constant weighting factor. Note, setting $\alpha = 0$ or $\alpha = 1$ introduces a non-trivial null space at frequencies that correspond to the interior resonance of the structure; the reader is referred to \cite{olaf_book} for a theoretical explanation. Here, $\mathcal{S}:H^{-1/2} (\Gamma) \longrightarrow H^{1/2} (\Gamma)$ denotes the single layer boundary integral operator
\begin{subequations}
	\begin{equation}
	\begin{aligned}
	\mathcal{S}[\Lambda,\Gamma] (\textbf{r}) &\doteq \int_{\Gamma}  \Lambda(\mathbf{r^\prime}) G(\mathbf{r},\mathbf{r^\prime})d\mathbf{r}^\prime \mbox{,     } \textbf{r} \in \Gamma ,
	\end{aligned}
	\label{eq:single}
	\end{equation}
	and $\mathcal{D}':H^{-1/2} (\Gamma) \longrightarrow H^{-1/2} (\Gamma)$ denotes the adjoint double layer operator 
	\begin{equation}
	\begin{aligned}
	\mathcal{D}'[\Lambda,\Gamma](\textbf{r})  &\doteq \left[\int_{\Gamma} \Lambda(\mathbf{r^\prime})
	\frac{\partial G(\mathbf{r},\mathbf{r^\prime})}{\partial \textbf{n}}d\mathbf{r}^\prime\right] \mbox{,      } \textbf{r} \in \Gamma.
	\end{aligned}
	\label{eq:hyper}
	\end{equation}
\end{subequations}
where $G(\mathbf{r},\mathbf{r'})=e^{-i\kappa|\mathbf{r}-\mathbf{r}'|}/4\pi|\mathbf{r}-\mathbf{r}'|$
is the free-space Helmholtz kernel in $\mathbb{R}^{3}$, $\kappa$ is the wavenumber, and $\beta = {i}/{\kappa}$. An $e^{i\omega t}$ dependence is assumed and suppressed. We introduce the far-field operator $\mathcal{L}_{far}:H^{-1/2} (\Gamma) \longrightarrow L^{2} (\Gamma_{s})$ given by
\begin{subequations}
	\begin{equation}
	\begin{aligned}
	\Phi^{s}_{m}(\hat{\textbf{r}}) \doteq \mathcal{L}_{far}[\Lambda,\Gamma] (\hat{\textbf{r}}) ,
	\end{aligned}
	\label{eq:farfield_system}
	\end{equation}
	\begin{equation}
	\begin{aligned}
	\mathcal{L}_{far}[\Lambda,\Gamma] (\hat{\textbf{r}}) \doteq  \frac{1}{4\pi} \int_{\Gamma} \Lambda(\mathbf{r^\prime}) e^{i\kappa\hat{\textbf{r}} \cdot \mathbf{r^\prime}}  d\mathbf{r}^\prime \mbox{,       } \hat{\textbf{r}} \in \Gamma_{s}.
	\end{aligned}
	\label{eq:farfield}
	\end{equation}
\end{subequations}
The exact solution of (\ref{eq:BurtonMiller}) is generally not available. The numerical solution of the integral equations is effected in a discrete setting using an isogeometric method based on Loop subdivision.

\section{Loop Subdivision based Isogeometric Method}

Next, we introduce Loop subdivision based isogeometric method that constructs subdivision surfaces for describing the geometry of the computational domain as well as representing the solution space for dependent variables. In this paper, we limit ourselves to the review of the elementary properties of subdivision surfaces; information provided in this section is purely for completeness and omits details that can be found in several Ref.~\cite{Ciraksubd,loop1987smooth,stam1998evaluation,JieDaultShankerChapterSubd,abdel_acoustics}.

\subsection{Subdivision Surfaces}
Let $\mathcal{T}^{k}$ denote a $k$-th refined control mesh, with vertices $V^{k} := \{\textbf{v}_{i}, i = 1,\ldots, N_{v}\}$ and triangular faces $P^{k} := \{\textbf{p}_{i}, i = 1,\ldots, N_{f}\}$; by definition, $\mathcal{T}^{0}$ denotes the initial control mesh. Without going into the details, we can represent a $C^{2}$ (almost everywhere) smooth limit surface  $\Gamma$, through an infinite number of iterative refinements of the control mesh following the loop subdivision scheme \cite{subdivision}. In practice, this prescription is not followed, i.e.,  there are closed form expressions for computing the limit surface $\Gamma$ for a given control mesh $\mathcal{T}^{k}$ in terms of quantities defined on the given control mesh \cite{stam1998evaluation}. 

\begin{figure}[!htpb]
	\centering
	\includegraphics[width=3cm]{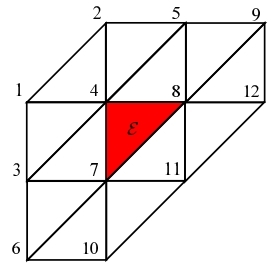}
	\caption{Regular triangular patch defined by its 1-ring vertices.
		\label{fig:1ring}}  
\end{figure}

Assume that a subdivision surface admits a natural parameterization of the surface $\Gamma$ in terms of the barycentric $(u,v)$ coordinates defined on each face $\epsilon \in P^{k}$, for some k. We begin by considering any patch $\epsilon \in P^{k}$ for some k, as depicted in Fig.~\ref{fig:1ring}. We define the 0-ring of a patch (triangle) as the vertices that belong to the patch, and the 1-ring as the set of all vertices, $n_{v}$, that can be reached by traversing no more than two edges, as shown in Fig.~\ref{fig:1ring}. We define the regularity of the triangle by the characterization of its vertices' valence (0-ring); the valence of a given vertex is the number of edges incident on itself. A vertex is considered regular if its valence is equal to 6, otherwise, it is called an irregular or extraordinary vertex. A triangle is regular if its vertices are all regular, and irregular otherwise. Using these definition, we can define the mapping from the barycentric coordinates $(u,v)$ that parameterize a patch $\epsilon$, to the limit surface $\Gamma(\textbf{r}(u,v))$ as 
\begin{equation}
\label{eq:Sexp}
\Gamma(\textbf{r}(u,v)) = \sum_{i=1}^{n_v} \textbf{c}_{i}\psi_{i}(u,v).
\end{equation}
where $\mathbf{c}_i$ are vertex locations of the $n_v$ control points. If a triangle is regular, then $n_{v} = 12$ and $\psi_{i}(u,v)$ is a box-spline basis function defined over the patch. Otherwise, when a triangle is irregular it must be refined until the considered point lies in a regular patch such that we can reapply the method described above \cite{stam1998evaluation}. Accordingly, we redefine $\psi_i$ as a subdivision basis set, and $N_{v}$ as the total number of control nodes defined on the 1-ring of the triangular patch.
In what follows, we have assumed that the mapping $\vb{r}(u,v)$ exists, where $(u,v)$ are the barycentric coordinates of a triangle; $\vb{r}$ will be used and the dependence on $(u,v)$ is assumed and suppressed.

\subsection{Properties of Loop subdivision basis\label{sec:sunbdivisionProperties}}
The Loop basis function's properties are critical to both the isogeometric analysis (IGA) as well as developing MHB. We enumerate them as such:

\begin{enumerate}
	\item Positivity: The basis functions associated with a control vertex are positive in its entire domain of support. 
	\item Compact support: The support of a basis function associated with a control vertex has compact support. For the scheme described here, this is the 1-ring associated with the vertex. Both, $\psi_l (\mathbf{r})$ and $\grad_s \psi_l (\mathbf{r})$ go smoothly to zero on the boundary. 
	\item Partition of Unity: The basis function forms a partition of unity. This implies that overlapping basis add to 1. 
	\item Continuity: The basis sets are $C^2$ everywhere except at irregular vertices where they are $C^1$. 
	\item Approximation theorem \cite{Arden}:
	\begin{theorem}
		Let $\mathcal{S}\left(\mathcal{T}^{k}\right)$ be the space of Loop’s subdivision functions on the k-times control mesh $\mathcal{T}^{k}$. For integers $0 \leq s < r \le 3$ and any $\epsilon > 0$ we have the following bound on the minimal $H^{s}\left(\mathcal{T}\right)$-approximation error of a function $f \in H^3 \left(\mathcal{T}\right)$:
		\begin{equation}
		dist\left(f,\mathcal{S}\left(\mathcal{T}^{k}\right)\right)_{H^{s}\left(\mathcal{T}\right)} \le C_{\epsilon}\lambda_{max}^{r-s-\epsilon}\norm{f}_{H^{r}\left(\mathcal{T}\right)},
		\end{equation}
		where the constant $C_{\epsilon} = C\left(\epsilon, \mathcal{T}\right)$ is independent of k and f and let $\lambda_{max} = \lambda_{max}\left(\mathcal{T}\right)$ be the largest subdominant eigenvalue of a subdivision matrix whose valence is represented by a vertex of $\mathcal{T}$ or 1/2 whichever is greater
	\end{theorem}
\end{enumerate}

\subsection{Loop subdivision based IGA }
Let $\Psi = \left \{ \psi_1 (\mathbf{r}), \psi_2 (\mathbf{r}), \cdots, \psi_{N_v}(\mathbf{r}) \right \} $ denote a set of Loop subdivision basis sets on $\Gamma$. It is apparent that $\Psi \subset H^{2}\left(\Gamma\right)$. The set of basis functions corresponding to each subdivided computable control patch are utilized to represent the geometry of interest and the solution space for dependent variables. Futhermore, these constructed functions seamlessly carry over the aforementioned properties into the physical domain. To this end, instead of control vertices, assume that there exists a net of control pressure (coincident with the location of the control net). The pressure, $\Phi (\mathbf{r})$, on the limit surface can then be expressed in terms of these basis sets via 
\begin{equation}
\label{eq:Hexp}
\Phi(\mathbf{r}(u,v))= \sum_{i=1}^{N_v} a_{i}\psi_{i}(\mathbf{r}(u,v)).
\end{equation}
where $a_{i}$ are  weights assigned to the locally indexed $i^{th}$ control vertex; $N_{v}$, $\psi_{i}(\mathbf{r}(u,v))$, and $(u,v)$ retain the same definition as those prescribed above. It is apparent that the set of basis functions that span $\Psi$ can be used to discretize the operators in \eqref{eq:BurtonMiller}. Using a Galerkin approach, we arrive at a system of equations $\mathcal{Z} \mathcal{I} = \mathcal{V}$, where
\begin{equation}
\begin{split}
\mathcal{Z}_{mn}  = \langle \psi_{m} (\mathbf{r}), \mathcal{L}_{BM}\left[\psi_{n},\Gamma\right](\mathbf{r}) \rangle \mbox{,     } \textbf{r} \in \Gamma, \\
\mathcal{V}  = \left[ v_1, v_2, \cdots, v_{N_{b}} \right ]^{T},    \\
\mathcal{I}  = \left[ a_1, a_2, \cdots, a_{N_{b}} \right ]^{T}, 
\end{split}
\label{eq:discBM}
\end{equation}
Here, $v_i = \langle \Psi_{m} (\mathbf{r}), \mathcal{V}^i(\mathbf{r}) \rangle = \int_{\Gamma} \psi_m (\mathbf{r}) \mathcal{V}^i(\mathbf{r}) d\textbf{r}$. Considering that the matrix $\mathcal{Z}_{mn}$ is dense an acceleration technique, such as wideband fast multipole method \cite{vikram2009novel}, are crucial for the fast and efficient computation of matrix operations. A detailed description of our implementation of IGABEM accelerated with multi-level fast multipole algorithm with a detailed error analysis is given in \cite{abdel_acoustics}.

\section{Manifold Harmonics using Loop Subdivision \label{sec:manifoldHarm}}

While we have formulated a Loop subdivison IGA framework, the costs are still prohibitively high when using it for shape reconstruction/optimization. As alluded to earlier, one workaround is utilizing the multi-resolution nature of subdivision \cite{ciraksubdsurfaces,Jan_multiresolution,FEM_multires,shape_FEM_multi}. The subdivision surfaces allow one to use different resolutions of the same geometry for optimisation and analysis by employing the hierarchy of the control mesh underlying a subdivision surface. In this paradigm, the degrees of freedom in optimization (i.e., design variables) are chosen as the vertex coordinates of a coarser control mesh and analysis is performed on the finer mesh. The mesh is evolved during optimization. This approach has been demonstrated for static problems, i.e., $\kappa = 0$ \cite{cirak_subd_statics}. While this approach is effective, it is apparent that this approach is perturbative and suffers from the same challenges as elucidated earlier; i.e., for geometrically complex objects, the number of degrees of freedom can rapidly increase leading to a computationally challenging problem.  

As a result, a more robust compression scheme that can be systematically enriched (in both geometry and physics) is necessary to make function evaluations more tractable. To achieve this, we employ MHB defined via eigenfunctions of the Laplace-Beltrami operator on $\Gamma$  by leveraging the representation power of Loop subdivison. As is well known, MHB constitutes a compact and stable basis for shape representation. Compactness means that most natural functions on a shape should be well approximated by using a small number of basis elements, while stability implies that the space of functions spanned by the basis functions must be stable under small shape deformations. These two properties together ensure that we can represent a function using a small and robust subset of MHs, implying we need only consider a subset of coefficients (i.e. set of weights assigned to each MH). The number of coefficients (i.e., design variables) is far fewer than the number of vertices leading to a highly compressed reconstruction/optimization scheme resulting in substantial gains in computational time. In what follows, we will discuss MHB in the context of Loop subdivision. 

\subsection{Laplace-Beltrami Operator}

To begin, let $H$ be a real-valued function defined on a compact 2D Riemannian manifold $\Gamma$ embedded in $\mathbb{R}^{3}$. The Laplace-Beltrami operator $\Delta$ is defined by
\begin{equation}
\Delta H := \nabla \cdot (\nabla H)
\end{equation}
The Laplacian eigenvalue problem is as follows: 
\begin{equation}
\label{eq:LBOequation}
\Delta H = -\lambda H
\end{equation}
Since the LBO is Hermitian, the eigenvalue spectrum of $\Delta$ acting on $\mathcal{L}_{2}(\Gamma)$ is  a countable set of nonnegative numbers $0=\lambda_{0}\le\lambda_{1}\le \ldots$  and $H_{1},H_{2},\ldots$ are the corresponding orthonormal eigenfunctions satisfying $\langle H_{i},H_{j} \rangle_{L_{2}(\mathcal{M})} = \delta_{ij} $; for a more detailed discussion on the main properties of the Laplace-Beltrami operator, we refer the reader to \cite{LBO_MH}.

\subsection{Subdivision FEM for Computing Eigenfunctions of LBO}

In order to numerically compute the eigenfunctions of the LBO or MHs, we discretize (\ref{eq:LBOequation}) using the Loop Subdivision FEM Galerkin method. This is akin to similar efforts using Lagrangian surface descriptions \cite{DLBO,LBO_DNA} that have shown both $h-$ and $p-$ convergence \cite{FEMstrang,DLBO,Fu2017GeneralizedDS,LBO_DNA}.  

The numerics necessary for computing eigenfuctions of the LBO relies on casting the Laplacian eigenvalue problem into a variational setting. The solution of this variational problem is approximated using the finite element Galerkin technique on the surface. We begin by multiplying Eq.(\ref{eq:LBOequation}) with some test function $v (\mathbf{r})$ and then use Green's theorems to arrive to the following:
\begin{equation}
\label{eq:weak_LBO}
\int_\Gamma  \left < \nabla_s v (\mathbf{r}), \nabla_s H(\mathbf{r}) \right > d\mathbf{r} = -\lambda \int_\Gamma  v(\mathbf{r}) H (\mathbf{r}) d \mathbf{r}.
\end{equation}
For the numerical computation of the Laplacian eigenvalues and eigenfunctions, a discretization of $H (\mathbf{r}) \in H_0^1 (\Gamma)$, as $H (\mathbf{r}) = \sum_i^{N_v} h_i \psi_i (\mathbf{r})$ for $h_i \in \mathbb{R}$ and choosing $v (\mathbf{r}) \in \Psi$ as the test function leads to the following the generalized eigenvalue problem 
\begin{equation}
\mathcal{A}\mathcal{H} = -\boldsymbol\Lambda \mathcal{B}\mathcal{H},
\label{eq:geneigvalprob}
\end{equation}
where,
\begin{subequations}
	\begin{equation}
	\begin{aligned}
	\mathcal{A}_{ij} = \int_\Gamma \nabla_s \psi_{i} (\mathbf{r}) \cdot \nabla_s \psi_{j} (\mathbf{r}) d \mathbf{r},
	\end{aligned}
	\label{eq:stiffmatrix}
	\end{equation}
	\begin{equation}
	\begin{aligned}
	\mathcal{B}_{ij} = \int_\Gamma \psi_{i} (\mathbf{r})  \psi_{j} (\mathbf{r}) d\mathbf{r}.
	\end{aligned}
	\label{eq:massmatrix}
	\end{equation}
\end{subequations}
For this generalized symmetric eigenvalue problem  $\mathcal{A} \in \mathbb{R}^{N_v \times N_v}$ is positive semi-definite, $\mathcal{B} \in \mathbb{R}^{N_v \times N_v}$ is positive definite, $\boldsymbol\Lambda \in \mathbb{R}^{N_v \times N_v}$ contains  $N_v$ eigenvalues along its diagonal, and $\mathcal{H} \in \mathbb{R}^{N_v \times N_v}$ is a column space, $\mathcal{H}_{1},\ldots,\mathcal{H}_{N_v}$, containing the coefficients of the  Laplacian eigenvectors; here each eigenvector is defined using $\mathcal{H}_i = \left [ h_{1,i}, \ldots, h_{N_v,i} \right]^T$. For this symmetric generalized eigenvalue problem we  have $\mathcal{H}^{T}\boldsymbol\Lambda\mathcal{H} = \boldsymbol\Lambda$ and $\mathcal{H}^{T}\mathcal{B}\mathcal{H} = \mathcal{I}$, where $\mathcal{I}$ is the idempotent. From the previous relations, it follows that the eigenfunctions of the geometric and FEM Laplacian matrix are orthogonal with respect to the $\mathcal{B}$-based scalar product (i.e., $\langle \mathcal{H}_{i}, \mathcal{H}_{j} \rangle_{L_{2}(\mathcal{M})}  = \mathcal{H}_{i}^{T}\mathcal{B}\mathcal{H}_{j}$). Purely for the sake of completion, it follows that there exist $N_v$ eigenfunctions defined by $H_m (\mathbf{r}) = \sum_i^{N_v} h_{m,i}\psi_i(\mathbf{r})$. The eigenvectors with corresponding eigenvalues can then be calculated with a direct eigensolver or by using the efficient band-by-band computation method presented in \cite{MHs}. There is a extensive body of literature on efficient computation of these functions, largely applied to computational graphics \cite{LBO_Zhang}.

\subsection{Analysis and Representation on/of Manifolds}

Thus far, we have defined a space of MHs $\left\{ \mathcal{H} = \left\{H_i\right\}_{i=1}^{N_{v}} : H_i \in H^2 (\Gamma) \right\}$ that inherit the properties of the subdivision basis sets. In what follows, we illustrate three critical features of this representation that will be useful for shape reconstruction. The first is the compressed representation of the manifold; the second, the ease with which the manifold can be manipulated; and the third, representation of functions on the manifold.

\subsubsection{Compressed Representation of the manifold} \label{Compressed Representation of the manifold}

Let $\Gamma$ be a 2D manifold embedded in $\mathbb{R}^{3}$. Considering $\mathcal{H}$, we seek a representation
\begin{equation}
\label{eq:MHT}
\Gamma (\mathbf{r}) = \sum_{i=1}^{N_v} \beta_{i}H_{i}(\mathbf{r}).
\end{equation}
It is apparent from the orthogonal property of MHs that the coefficients can be defined as $\beta_{i} = \langle \Gamma (\mathbf{r}) ,H_{i} (\mathbf{r}) \rangle_{L_{2}(\mathcal{M})} = \int_\Gamma d \mathbf{r} \Gamma (\mathbf{r}) H_i (\mathbf{r})$; this, in effect, is dubbed a Manifold Harmonic Transform (MHT) of the surface, see \cite{MHs}. Furthermore, it is important to remember that each eigenfunction $H_i (\mathbf{r})$ corresponds to different spatial ``frequency'' on the manifold. Thus, the magnitude of $\beta_i$, dictate the relative importance of a given eigenfunction at a given spatial ``frequency''. This transformation makes it possible to define equivalent signal processing functions; these functions include, concepts such as under/over sampling, filtering, multi-resolution, windowing, and so on.  To wit, this representation makes it  possible to enrich shapes (and functions defined on the shape) systematically. To illustrate, the effectiveness of this approach, consider the manifold $\Gamma (\mathbf{r})$ again. Equation \eqref{eq:MHT} can be rewritten as 
\begin{equation}
\label{eq:manifold_coeffs}
\Gamma (\mathbf{r})  = \underbrace{\sum_{i=1}^{M}\beta_{i}H_{i}(\mathbf{r})}_\text{Low Frequency} 
+ \underbrace{\sum_{i=M+1}^{N_v}\beta_{i}H_{i}(\mathbf{r})}_\text{High Frequency}. 
\end{equation}
In this expression, we have designated some cutoff coefficient $\beta_M$, or correspondingly, $\sqrt{\lambda_M}$ to be a maximum cutoff ``frequency''. This effectively is a versatile tool for shape reconstruction that enables us to enrich data at different levels of fidelity. 

\begin{figure}[!h]
	\centering
	\vspace{1cm}
	\subfloat[Original model.]{
		\label{girl_original_modes}
		\includegraphics[width=0.49\textwidth]{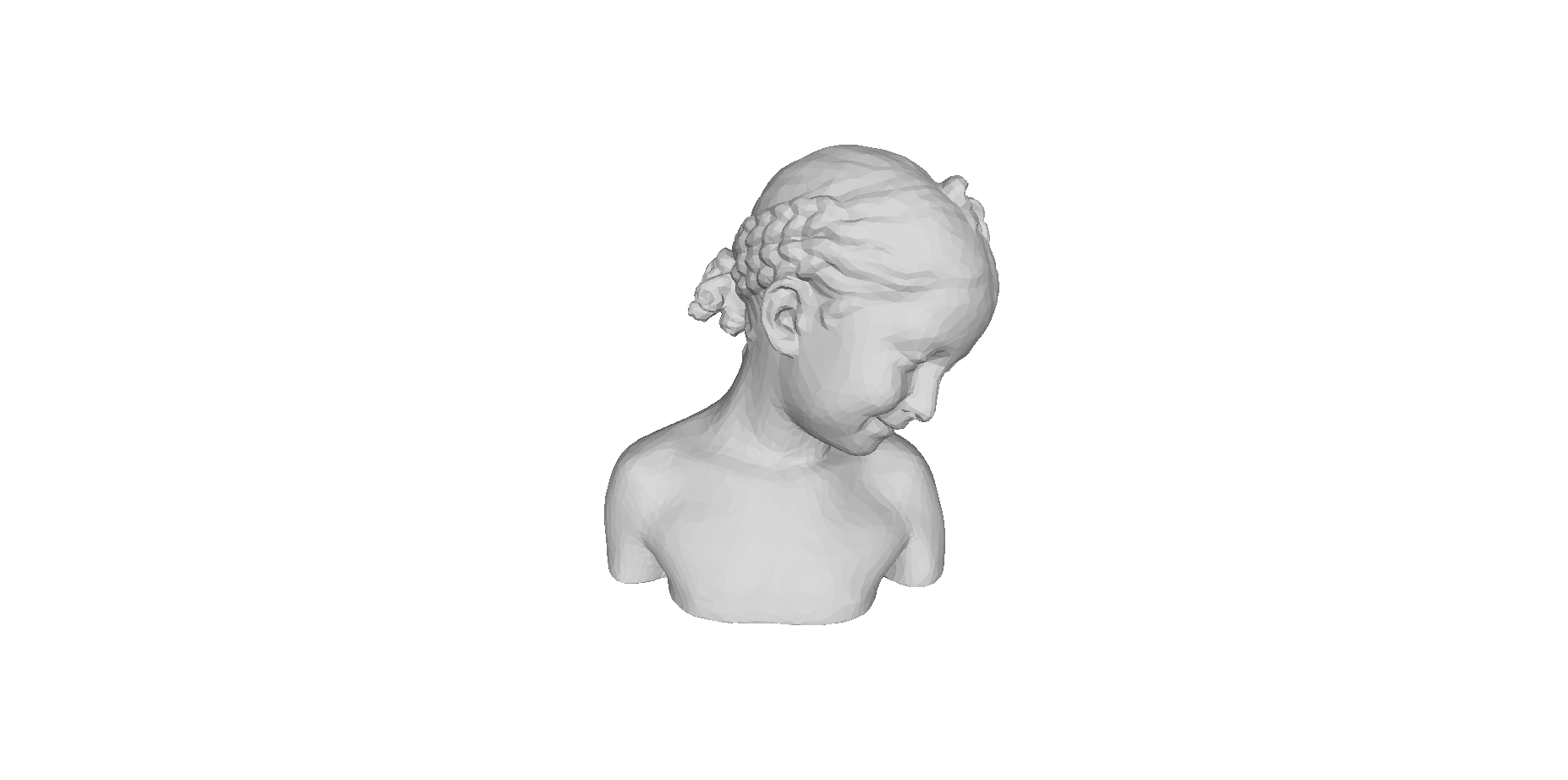} } 
	%
	\subfloat[With 100 MHs.]{
		\label{girl_100_MHS}
		\includegraphics[width=0.49\textwidth]{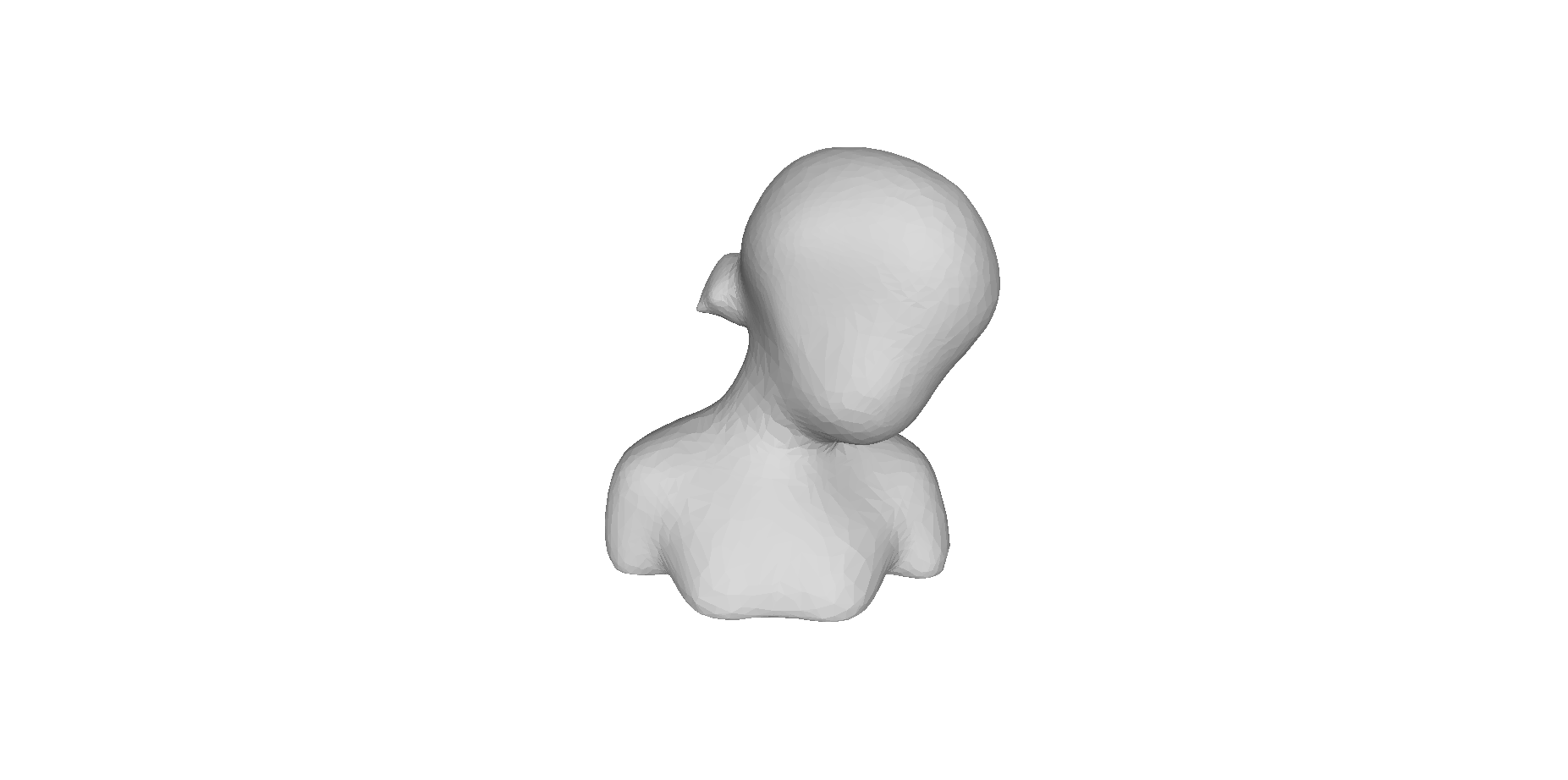}} 
	
	\subfloat[With 250 MHs.]{
		\label{girl_250_MHS}
		\includegraphics[width=0.49\textwidth]{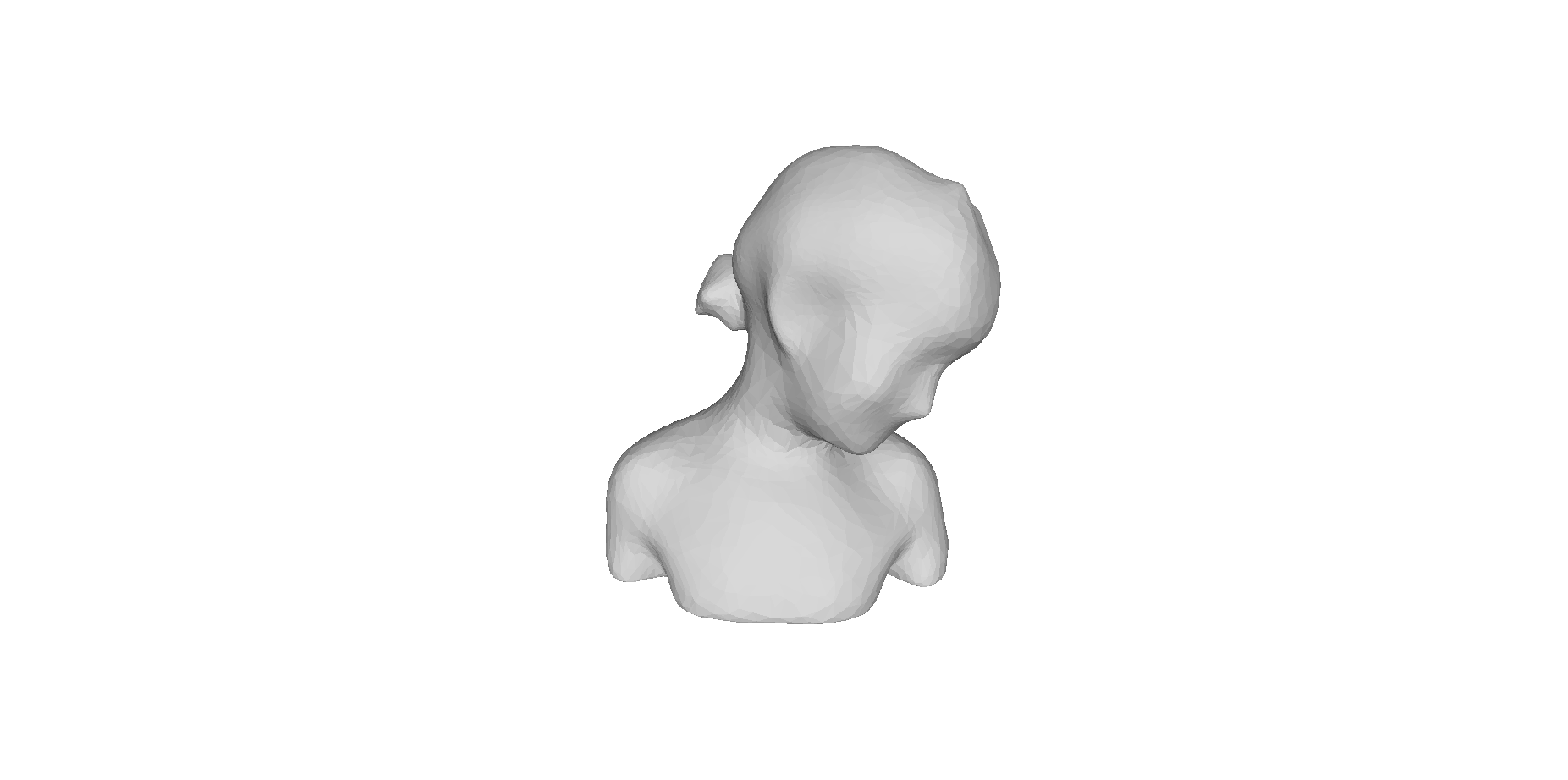}} 
	%
	\subfloat[With 750 MHs.]{
		\label{girl_750_MHS}
		\includegraphics[width=0.49\textwidth]{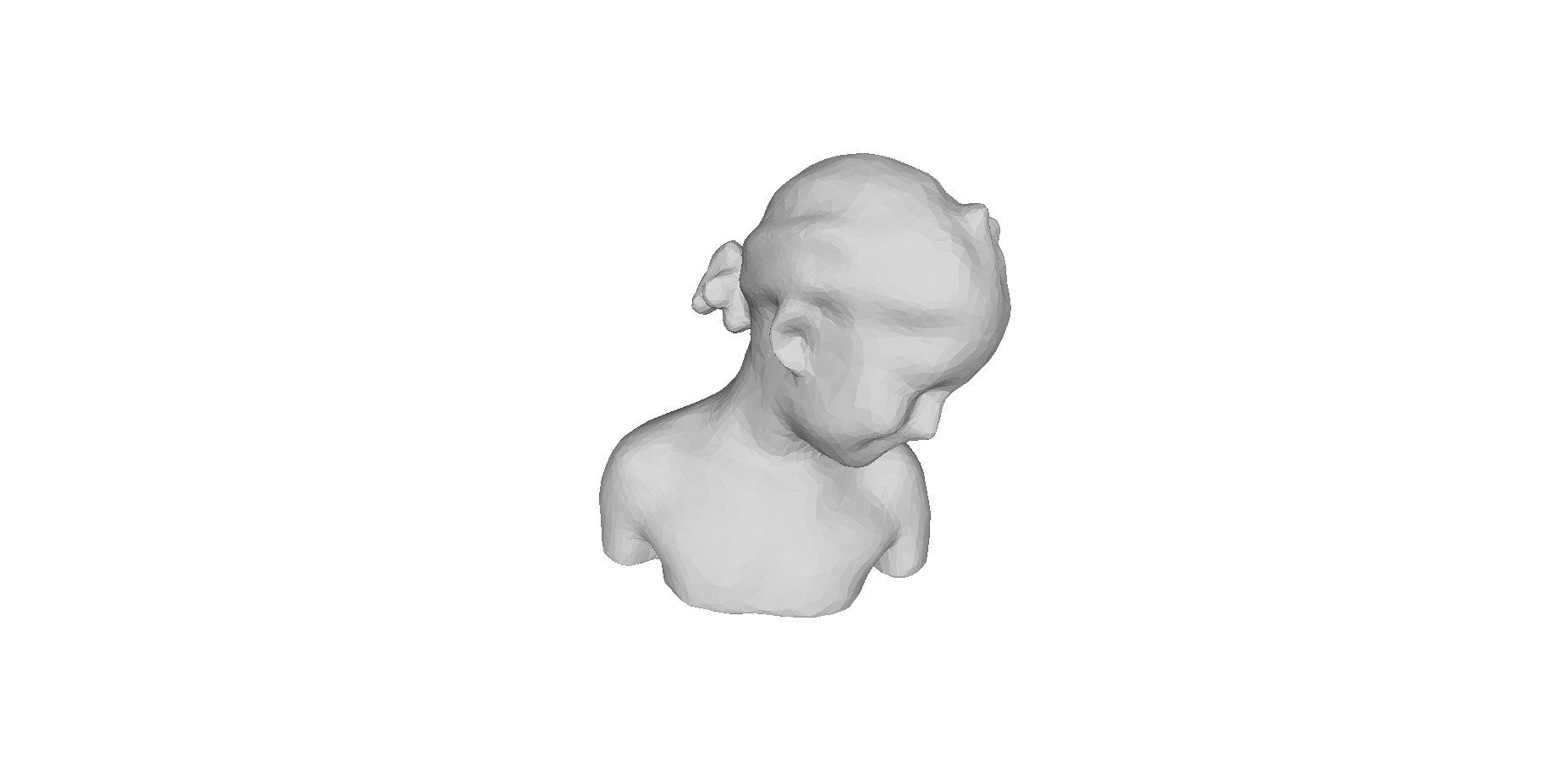} } 
	\caption{The statue of the girl in (\ref{girl_original_modes}) is reconstructed with increasing number of MHs in (\ref{girl_100_MHS})-(\ref{girl_750_MHS})}
	\label{fig:girl_MH_reconstruction}
\end{figure}

Figure \ref{fig:girl_MH_reconstruction} depicts the reconstruction of a geometry as we increase the number of MHs. As is evident from these figures, it does not take too many MHs to capture the general shape of the object. Furthermore, as we increase the number of MHs, we rapidly approach the true shape, i.e. we add more localized features. Note, the total number of MHs in the original system is $N_v = 5002$.

In this paper, we aim to utilize the natural multi-resolution framework of MHB detailed above for the shape reconstruction problem. Specifically, lets denote the set of coefficients in \eqref{eq:manifold_coeffs} using $\left \{ \beta \right \} = \left[\beta_1, \cdots, \beta_{N_v} \right ] $. It is apparent that one can group these in a collection of spatial frequency bands. In effect, $\left \{ \beta \right \} = \left [ \mathcal{B}_1, \cdots, \mathcal{B}_K \right ]$ where $\mathcal{B}_i  = [\beta_{i,1}, \cdots, \beta_{i,n_i} ]$, where there are $n_i$  coefficients in each set. Note, the sets are a contiguous partition of $\left \{ \beta \right \}$. It follows, that one can implement a multi-resolution analysis/optimization scheme by updating collection of $\mathcal{B}_i$ one at a time, as opposed to the entire set.

\subsubsection{Analysis using manifold harmonics}
The MHs $\mathcal{H}$ defined earlier can be utilized as an isogeomtric basis set enabling us to not only represent  geometry, but as a basis for representing physical quantities defined on the geometry; we note that the optimilaity of MHs for signal approximation has been addressed \cite{LBO_signal_rep}. The salient features of this representation are as follows: (a) Subdivision-based MHs are $C^{2}$ smooth allowing for an excellent basis for representing physics on the geometry; (b) the orthogonal property allows for hierarchical signal representation, de-noising, and compression; (c) low computational cost and storage overhead feasible for applications in analysis; and (d) subdivision-based MHs span the same space as $\Psi$. 

Considering the aforementioned properties, we aim to construct a spectral representation of scalar functions defined on the surface using MHs. Following Eq.~\ref{eq:MHT}, we have $\Phi(\mathbf{r}) = \sum_{i=1}^{N_v} \alpha_{i}H_{i}(\mathbf{r})$, wherein it is possible to choose an error threshold, $\epsilon$, for representation of desired functions on the manifold and thereby, truncate the number of MHs used. To asses the efficiency and accuracy of our MH spectral representation method for both physics and geometry we consider a $3.21\lambda \times 3.21\lambda \times 3.04\lambda $ bumpy cube in Fig.~\ref{fig:bumpy_cube_exact} as our candidate object; this object is represented using 2562 subdivision basis functions. The geometry is represented using subdivision basis. The physics on this surface is represented using an increasing number of MHs, and solution obtained is compared against an isogeometric solution for a plane wave incident along the $\hat{z}$ direction. We define the error between the full isogeometric solution to those using MHs 
\begin{equation}
\epsilon_{L_{2}} = \frac{ \left\Vert   \tilde{\Phi} - \Phi \right\Vert}{\left\Vert \Phi \right\Vert}
\end{equation}
where $\left\Vert  \cdot \right\Vert$ denotes the $L_{2}$ norm, $\Phi$ denotes the scattered far field obtained using purely the subdivision basis for physics and  $\tilde{\Phi}$ the solution obtained using MHs. As is evident from Table.~\ref{table1}, we  observe convergence with increasing MHs. 

\begin{figure}[!h]
	\centering
	\subfloat[10 Geometric MHs]{
		\label{fig:bumpy_cube_10}
		\includegraphics[width=0.30\textwidth]{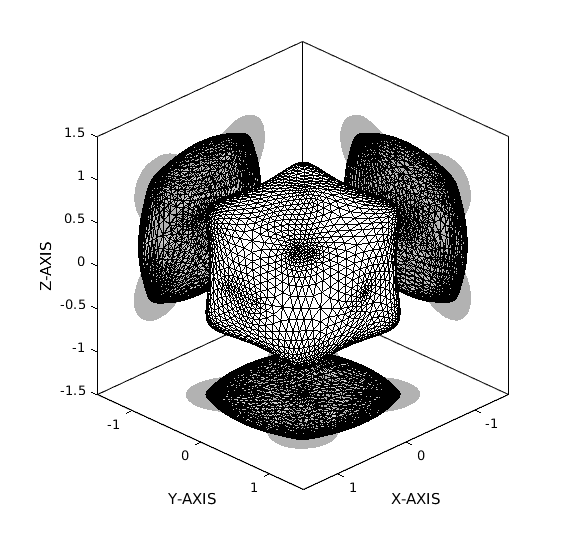} } 
	\hspace{1cm}
	\subfloat[50 Geometric MHs]{
		\label{fig:bumpy_cube_50}
		\includegraphics[width=0.30\textwidth]{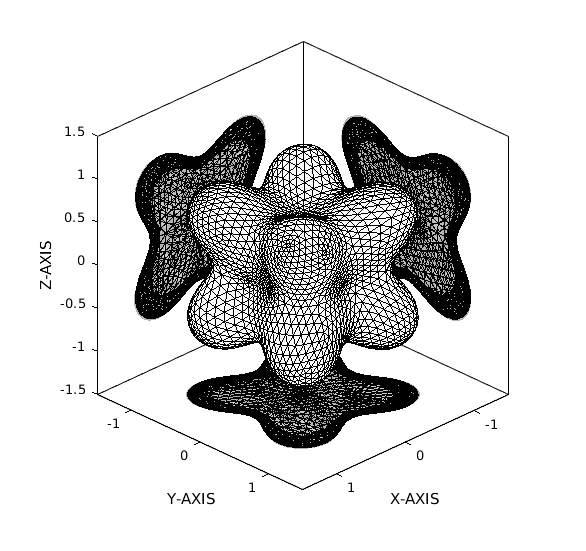} }
	\vspace{0.8cm}
	\\
	\subfloat[100 Geometric MHs]{
		\label{fig:bumpy_cube_100}
		\includegraphics[width=0.30\textwidth]{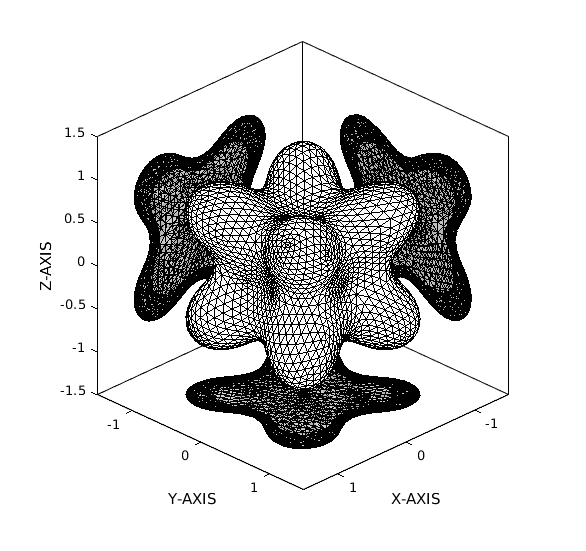}} 
	\hspace{1cm}
	\subfloat[200 Geometric MH]{
		\label{fig:bumpy_cube_200}
		\includegraphics[width=0.30\textwidth]{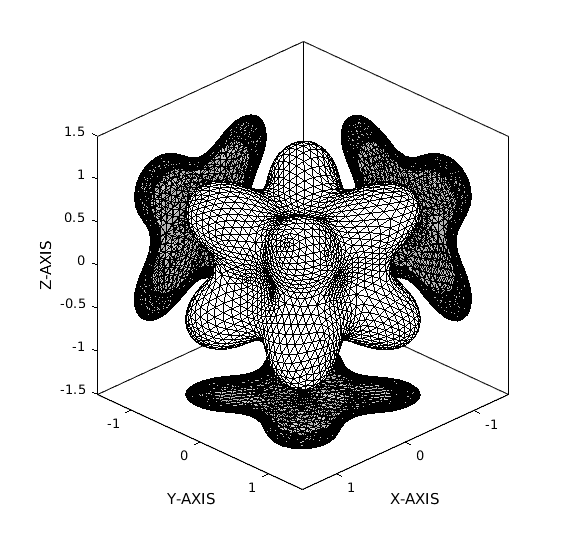} } 
	
	\subfloat[$\epsilon_{L_{2}}$ relative error]{
		\label{fig:phy_geom_MHs}
		\includegraphics[width=0.45\textwidth]{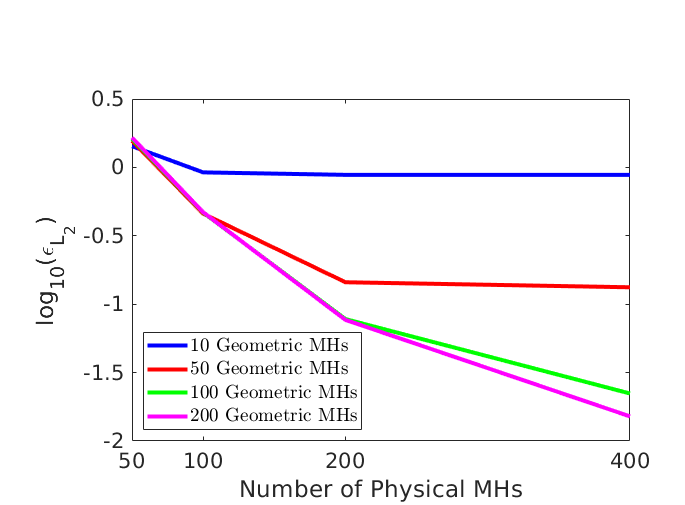} } 
	\caption{ Convergence in $\epsilon_{L_{2}}$ relative error of reconstructed far field as we increase MHs used to represent the surface and fields on the surface}
	\label{fig:analysis_MHs}
\end{figure}

Next, we study the reconstruction of both the geometry and the physics using MHs in tandem. In the first row in Table.~\ref{table_1b}, we present the geometric errors with respect to the number of MHs. The metrics we will use here and throughout the paper are as follows: 
(a) Surface area error $S_{err}\left(\Gamma_{G},\Gamma_{n}\right) = \frac{\left \Vert  S_{n}-S_{G} \right \Vert}{\left \Vert S_{true} \right \Vert}$, where $S_{n}$ is the reconstructed surface area for $n$ MHs and  $S_{E}$ is the exact surface area. And, (b) $H\left(\Gamma_{E},\Gamma_{n} \right)$, which denotes the maximum Hausdorff distance \cite{meshlab}, between the exact geometry and the reconstructed geometry for $n$ MHs. In addition, in Figs.\ref{fig:bumpy_cube_10}-\ref{fig:bumpy_cube_200} both the reconstructed surface as well as the reconstructed and exact surfaces are prejected onto the $x$-$z$, $y$-$z$, and $x$-$y$ planes for direct comparisons, indicated as black meshed and gray solid regions, respectively. 

Finally, Fig.\ref{fig:phy_geom_MHs} demonstrates convergence as one increases the number of MHs for both geometry and physics. Each of the plots corresponds to fixed MHs for geometry with increasing MHs for physics. The error is measured against those obtained via an isogeometric solve. As is evident, it is possible to represent both the geometry and scattered field using significantly fewer MHs. It should also be noted that the errors plateau soon after a threshold error is reached in geometry reconstruction. 

\begin{table}[!h]
	\begin{center}
		\begin{tabular}{  m{2cm} | m{1.5cm} | m{1.5cm} | m{1.5cm}| m{1.5cm} |m{1.5cm}  } 
			\hline
			No. of MHs & 500 & 1000 & 1500 & 2000 & 2562 \\ 
			\hline
			Rel. $\epsilon_{L_{2}}$ error & 6.19E-4 & 3.50E-5 & 1.58E-5 & 9.55E-6 & 1.25E-13 \\ 
			\hline
		\end{tabular}
	\end{center}
	\caption{$\epsilon_{L_{2}}$ relative error for spectral representation of scattered fields for exact geometry.}
	\label{table1}
\end{table}

\begin{table}[!h]
	\begin{center}
		\begin{tabular}{  m{2cm} | m{1.5cm} | m{1.5cm} | m{1.5cm}| m{1.5cm}  } 
			\hline
			No. of MHs & 10 & 50 & 100 & 200 \\ 
			\hline
			$S_{err}\left(\Gamma_{G},\Gamma_{n}\right)$ & 2.46E-1 & 1.89E-2 & 4.13E-3 & 2.32E-3  \\   
			\hline 
			$H\left(\Gamma_{G},\Gamma_{n} \right)$ & 3.10E-1 & 4.24E-2 & 3.97E-2 & 3.30E-2 \\ 
			\hline
		\end{tabular}
	\end{center}
	\caption{Geometric errors for spectral representation of geometry.}
	\label{table_1b}
\end{table}

\subsection{Back-projection Based Initialization}
\label{Inverse Sources as an initial guess}
Next, we need to obtain an initial guess for the shape  $\Gamma$ from a set of target far-fields available from a multitude of incident angles and frequencies. To utilize this data, we use well developed techniques in volume source reconstruction method to obtain a point cloud that can then be used to construct a subdivision mesh. Let us denote the available data as $\Phi^s \left (\kappa_n, \hat{\kappa}_m, \hat{\mathbf{r}} \right )$ where $\kappa_n = \omega_{n}/c$ is the wave number, $\hat{\kappa}_m$ is an incident wave direction and $\hat{\mathbf{r}}$ is the unit-vector along the direction of observation (alternatively defined in terms of $\theta$ and $\phi$). Following \cite{Alvarez2012TheSR}, we aim to estimate the boundary of the object, by assuming there exists some source distribution $\Lambda_{eq} \left (\kappa_n, \hat{\kappa}_m, \mathbf{r} \right )$  in a volume $V$ that satisfy the following minimization problem 
\begin{equation}
\substack{min \\ n,m} \parallel \Phi^{s} \left (\kappa_n, \hat{\kappa}_m, \hat{r} \right ) - \mathcal{L}_{far} \left[ \Lambda_{eq}\left ( \kappa_{n}, \hat{\kappa}_m, \mathbf{r} \right ), V \right ]  \parallel
\end{equation}
Here $\mathcal{L}_{far} \left [ \Lambda_{eq}, V\right ]$ is a far-field projector similar to \eqref{eq:farfield}; the domain of integration is defined over $\Gamma_{s}$ and the function is defined over the volume $V$. We then coherently combine these volume equivalent sources in a manner following \cite{SAR_coherent},
\begin{equation}
\tilde{\Phi}_{eq}(\textbf{r}) = \sum_{m=1}^{M_{inc}}\sum_{n=1}^{N_{freq}}\Phi^{m,n}_{eq}(\textbf{r})  e^{i\kappa_{n}(\textbf{r} \cdot \hat{\vb{k}}_{m})},\hspace{0.5cm} \vb{r} \in V.
\label{eq:coherent_sum}
\end{equation}

\subsubsection{Initial Guess example}
\begin{figure}[!h]
	\centering
	\subfloat[Original model.]{
		\label{fig:spot_org}
		\includegraphics[width=0.20\textwidth]{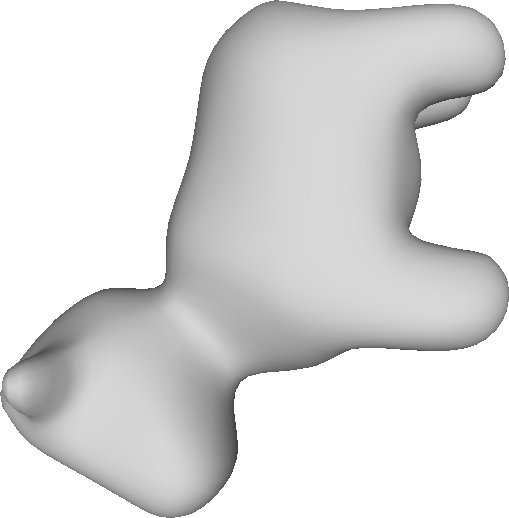} } 
	\hspace{2cm}
	\subfloat[Cross section of 3-D reconstructed object]{
		\label{fig:spot_2d}
		\includegraphics[width=0.21\textwidth]{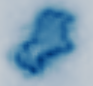} } 
	
	\subfloat[Point cloud]{
		\label{fig:spot_point_cloud}
		\includegraphics[width=0.20\textwidth]{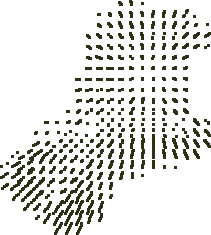}} 
	\hspace{2cm}
	\subfloat[Mesh generated from point cloud]{
		\label{fig:spot_cloud_to_mesh}
		\includegraphics[width=0.27\textwidth]{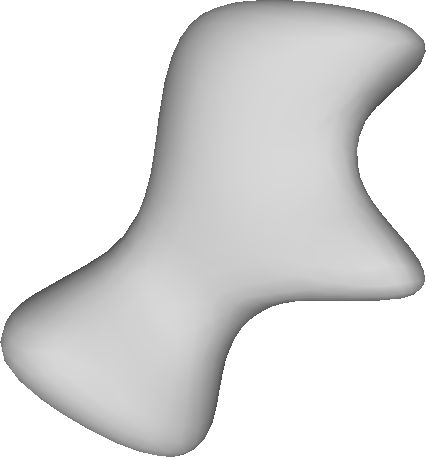} } 
	\caption{Initial guess pipeline: VSRM }
	\label{fig:MH_reconstruction}
\end{figure}

We will illustrate how to construct a starting BEM subdivision mesh using the initial guess routine detailed above for the shape reconstruction problem. Our candidate object is the cow in Fig.~\ref{fig:spot_org}; we are given a set of target farfield scattering patterns obeying Eq.~\ref{eq:forward_problem} due to 16 incident waves operating at a range of 200 Hz to 1200 Hz in 50 Hz increments. The reconstruction domain V is a cubic region of 4 m x 4 m x 4 m discretized into 50 x 50 x 50 grid points. At each point, we coherently sum the back propagated fields leading to an equivalent acoustics pressure in the volume. The points that are of the highest intensity are closest to the boundary of the object $\Gamma$; this effect can be seen in Fig.~\ref{fig:spot_2d}, wherein a cross-section of the reconstructed equivalent acoustic pressure along the $z$ axis is plotted. At this point we prescribe some threshold such that we isolate the points nearest to the boundary of the target object. This leads to a point cloud in the shape of the test object, see Fig.~\ref{fig:spot_point_cloud}, which is used to construct the initial mesh \ref{fig:spot_cloud_to_mesh}.

\subsection{Shape Reconstruction Algorithm}
Thus far, we have developed/described a comprehensive set of tools that we employ for shape optimization. The reconstruction proceeds as follows: given initial farfield data, we (a) obtain an initial mesh, (b) construct MHs for this initial guess, i.e. the intial set of coefficients $\beta_{0}$, and (c) use a optimization routine to find the set of coefficients that minimizes the constrained optimization problem in Eq.~\eqref{eq:objective_function} wherein  Eq.~\eqref{eq:forward_problem} is the constraint on the admissible fields. The steps for our shape reconstruction algorithm are outlined in algorithm  \ref{alg:shapeoptimzation}. In addition, we note that the optimization method used in this paper is the Method of Multiple Asymptotes (MMA) presented from the NLopt library \cite{NLopt}. However, the infrastructure presented in this paper is \emph{agnostic} to the specific optimization methodology used. 

\section{Numerical Results}
We present a series of tests that validate the approach presented in this paper. Our goal is to reconstruct complex simply connected shapes from synthetic data. As stated, we will assume that far field data is available for a set of incident illuminations and frequencies. An incident illumination, $\mathbf{d}_i$, is identified by the pair $(\theta_i, \phi_i)$. The frequency of the incident field with wavenumber $\kappa_j$ is denoted using $f_j$. In all cases presented, the solution outline is as follows: (a) we obtain synthetic data, or goal fields, from a subdivision description of the scatterer, (b) we add noise to this data defined by 
\begin{equation}
\hat{\Phi}^{s}_{G} (\mathbf{r}) = \Phi^{s}_{G}  (\mathbf{r}) + \delta \Phi^{s}_{G}  (\mathbf{r})
\end{equation}
where ${\hat{\Phi}}^{s}_{G}  (\mathbf{r})$ is the set of perturbed goal fields. The term $\delta >0$ is related to the signal-to-noise ratio (SNR) by the equation
\begin{equation}
\delta = \frac{10^{-SNR/20}}{\sqrt{M}},
\end{equation}
and is a unit vector with $M$ random Gaussian entries. The $SNR$ (in dB)  typically fluctuates in real applications between 10 dB--40 dB. For all our experiments, we choose a noise level of 10 dB $SNR$ to demonstrate the robustness of our algorithm to noise; a $SNR$ of 10 dB is considered considered a high level of noise. Once we have our new perturbed goal fields, we follow employ our shape reconstruction routine following the steps presented in Algorithm \ref{alg:shapeoptimzation}.

\begin{algorithm}[!h]
	\begin{algorithmic}[1]
		\caption{Shape reconstruction algorithm}
		\label{alg:shapeoptimzation}
		\STATE{\textbf{Define:}}
		\bindent
		\STATE{ $\mathbf{D}$: Directions of propagation}
		\STATE{$\vb{K}$: Wavenumbers}
		\STATE{$\Phi^i$: Incident fields for a set $\left\{\mathbf{D},\vb{K}\right\}$ }
		\STATE{$\Gamma$: The boundary to be determined}
		\STATE{$\Phi^{s}_{G}$: Goal fields at $\Gamma_{s}$ due to incident fields $\Phi^{s}$ impinging on $\Gamma$}
		\STATE{Termination tolerance $\epsilon$ for optimization routine}
		\eindent
		\STATE{\textbf{Method:}}
		\bindent
		\STATE{\textbf{VSRM: The Initial Guess}}
		\STATE{Given $\Phi^{s}_{G}$, construct the volumetric acoustics source $\Phi_{eq,V}$ using  Eq.~\ref{eq:coherent_sum}}
		\STATE{Generate an analysis-ready subdivision mesh $\mathcal{M}_{0}$, from the $\Phi_{eq,V}$ point cloud using \cite{pcloud_surf_recon}} 
		\STATE{ Solve Eq.~\ref{eq:geneigvalprob} on $\mathcal{M}_{0}$ to generate the MHs}
		\STATE{Construct $\beta_{0}$ by performing a MHT on $\mathcal{M}_{0}$}
		\STATE{\textbf{Optimization routine:}}
		\WHILE{$J(\beta_{i}) \geq \epsilon$}
		\bindent
		\STATE{Solve the forward problem (Eq. (1)–(3)) on $\mathcal{M}_{i}$}
		\STATE{Compute objective function $J(\beta)$ and corresponding sensitivity analysis}
		\STATE{Update design variables using an optimization algorithm: $\left \{ \beta \right \} _{i} \rightarrow \left \{ \beta \right \}_{i+1}$}
		\STATE{Update the subdivision mesh: $\mathcal{M}_{i} \rightarrow \mathcal{M}_{i+1}$}
		\eindent
		\ENDWHILE
		\eindent
	\end{algorithmic}
\end{algorithm} 

Note, every step in the shape reconstruction algorithm requires movement of the scatterer’s surface. Therefore, we must account for large deformations of the surface throughout the evolution of the algorithm. In order to do so, we require a high fidelity mesh. Given that we are using subdivision surface, we have $\approx \lambda / 7$ patches per wavelength. Furthermore, we note that in the case of a good initial guess, large mesh deformation will be far less prevalent and therefore easing the restriction on the fidelity of the mesh. Recall, that the goal is to achieve a reduction in the cost functional and not to find the global minimum of the non-convex optimization problem. Lastly, in all examples we provide both the surface area error,  $S_{err}\left(\Gamma_{G},\Gamma_{n}\right)$, the Hausdorff measure $H\left(\Gamma_{G},\Gamma_{n} \right)$ and the error in the functional. Furthermore, to illustrate the difference between the reconstructed and the goal surfaces we
added to the figures their projections onto the $x$-$z$, $y$-$z$, and $x$-$y$ planes, which are indicated as black meshed and gray solid regions, respectively.

\subsection{First Example: Bean}
\begin{figure}[!t]
	\centering
	\subfloat[]{
		\label{fig:bean_initial_guess}
		\includegraphics[width=4.5cm,height=3.5cm]{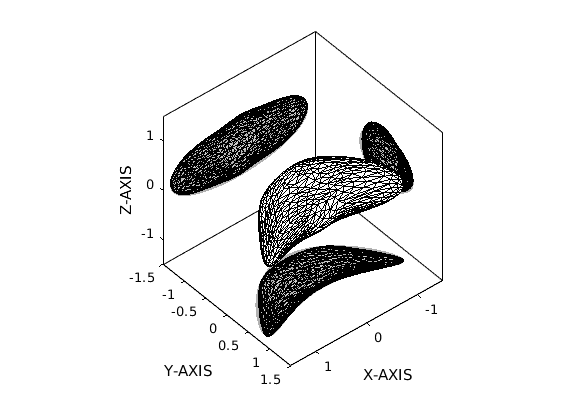} } 
	\hspace{1cm}
	\subfloat[]{
		\label{fig:sphere_inital_shape}
		\includegraphics[width=4.5cm,height=3.5cm]{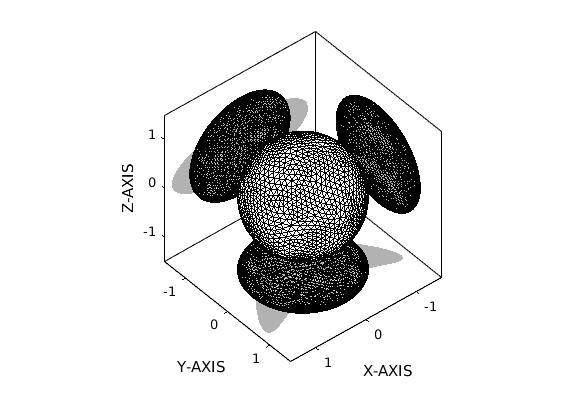} }
	\\
	\subfloat[]{
		\label{fig:VSRM_bean_final_shape}
		\includegraphics[width=4.5cm,height=4.5cm]{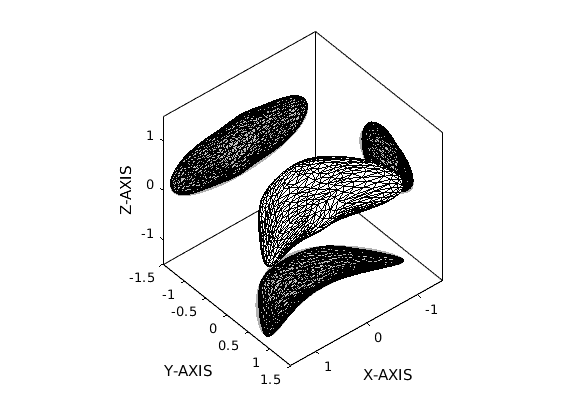}} 
	\hspace{1cm}
	\subfloat[]{
		\label{fig:SPH_bean_final_shape}
		\includegraphics[width=4.5cm,height=3.5cm]{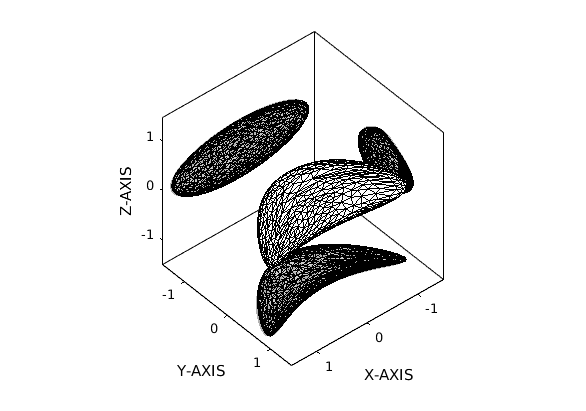} } 
	\\
	\subfloat[]{
		\label{fig:bean_exact_surface}
		\includegraphics[width=4.5cm,height=3.5cm]{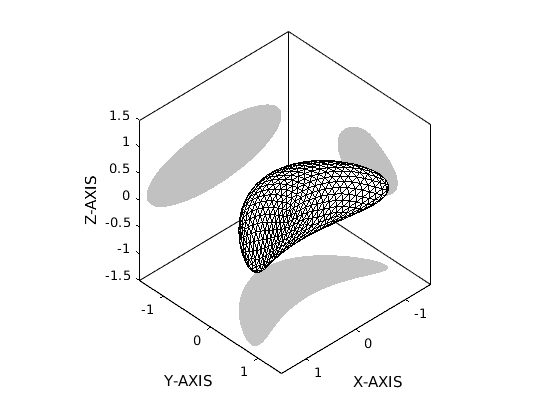}} 
	\caption{Shape reconstruction of a bean (\ref{fig:bean_exact_surface}): from an initial guess (\ref{fig:bean_initial_guess}), after 33 iterations we have (\ref{fig:VSRM_bean_final_shape}), $J(\beta) = 6.60 \cdot 10^{-3}$; from a sphere after 33 iterations (\ref{fig:SPH_bean_final_shape}) with $J(\beta) = 3.96\cdot 10^{-2}$.}
	
\end{figure}

In the first example, we consider a bean that fits in a  2.71 m $\times$ 1.27 m $\times$ 1.18 m bounding box as our target shape, see Fig.~\ref{fig:bean_exact_surface}. This target shape is sufficiently complicated in that it poses a challenging shape reconstruction problem, but still remains feasible for comparison of an arbitrary starting point against a VSRM initial guess. In particular, we conduct two different experiments: the first is to reconstruct the target shape starting from an arbitrary starting point, in this case a sphere of radius 1.0 m and the second, is to use VSRM to generate an initial guess, see Fig.~\ref{fig:bean_initial_guess}. When starting with a sphere, our goal fields are generated using the procedure described earlier using 9 incident fields at a set of directions $\mathbf{D}_{1} = \left \{ (\theta_1, \phi_1) \right \}$. The frequency of these fields are $f_1 = 540 Hz$. The initial starting mesh for the sphere is constructed using 1802 vertices and 3600 faces. The reconstruction is performed using 10 MHs.

The second experiment is to start from a VSRM guess; we use 16 incident plane waves defined by $\mathbf{D}_{2} = \left \{ (\theta_{2}, \phi_{2}) \right \}$, and frequencies $\left \{100 Hz, 125 Hz, 150 Hz, \cdots, 700 Hz\right \}$ to generate this initial guess. The starting mesh using the initial guess is discretized at 1202 vertices and 2400 faces. Similarly, the goal field used for the reconstruction process is generated from a single field operating at $f_2 = 540 Hz$ traveling in the $-\hat{\textbf{y}}$ direction. Lastly, we use 20 MHs for the reconstruction routine.

As seen in Figs.~\ref{fig:VSRM_bean_final_shape} and \ref{fig:SPH_bean_final_shape}, the algorithm was able to reconstruct the object and its features accurately from noisy data; notably the concavity of the object were successfully recovered in both cases. It should be noted that the algorithm was able to recover the correct size and width of the scattering object. Their relative placement and generic dimensions agree very well with those of the target.

\begin{table}[!h]
	\begin{center}
		\begin{tabular}{  m{3.0cm} | m{1.5cm} | m{1.5cm} | m{1.5cm}| m{1.5cm} |m{1.5cm}  } 
			\hline
			No. of iter. & Initial & 10 & 20 & 25 & 33 \\ 
			\hline
			$H\left(\Gamma_{G},\Gamma_{n} \right)$ & 5.37E-1 & 1.96E-2 & 8.21E-2 & 4.23E-2 & 4.87E-2 \\
			\hline
			$S_{err}\left(\Gamma_{G},\Gamma_{n}\right)$ & 5.75E-1 & 8.06E-2 & 3.15E-2 & 5.15E-3 & 8.97E-3 \\ 
			\hline
			$J(\beta_{n})$ & 6.08 & 6.04E-1 & 1.42E-2 & 9.15E-2 & 3.96E-3 \\ 
			\hline
		\end{tabular}
	\end{center}
	\caption{Evolution of the error metrics with respect to iteration for the bean Ref~(\ref{fig:bean_exact_surface}) starting from Ref~(\ref{fig:sphere_inital_shape}).}
	\label{table2}
\end{table}

\begin{table}[!h]
	\begin{center}
		\begin{tabular}{  m{3.0cm} | m{1.5cm} | m{1.5cm} | m{1.5cm}| m{1.5cm} |m{1.5cm}  } 
			\hline
			No. of iter. & Initial & 10 & 20 & 25 & 33 \\ 
			\hline
			$H\left(\Gamma_{G},\Gamma_{n} \right)$ & 1.91E-1 & 1.73E-1 & 8.26E-1 & 8.49E-2 & 8.03E-2 \\
			\hline
			$S_{err}\left(\Gamma_{G},\Gamma_{n}\right)$ & 1.24E-1 & 8.28E-2 & 3.45E-2 & 3.07E-2 & 2.28E-2 \\ 
			\hline
			$J(\beta_{n})$ & 2.57E-1 & 1.09E-1 & 1.06E-2 & 8.14E-2 & 6.60E-3 \\ 
			\hline
		\end{tabular}
	\end{center}
	\caption{Evolution of the error metrics with respect to iteration for the bean Ref~(\ref{fig:bean_exact_surface}) starting from Ref~(\ref{fig:bean_initial_guess}).}
	\label{table3}
\end{table}

\subsection{Second Example: Bumpy Cube}

\begin{figure}[!ht]
	\centering
	
	\subfloat[]{
		\label{fig:bumpy_cube_initial}
		\includegraphics[width=5cm,height=5cm]{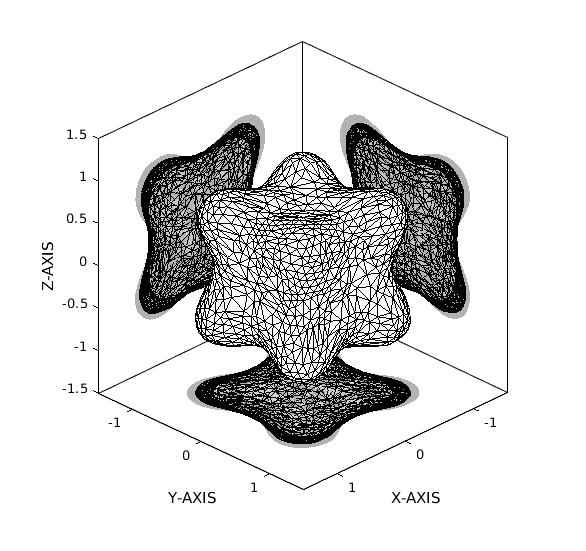} } 
	\hspace{1cm} 
	\subfloat[]{
		\label{fig:bumpy_cube_final}
		\includegraphics[width=5cm,height=5cm]{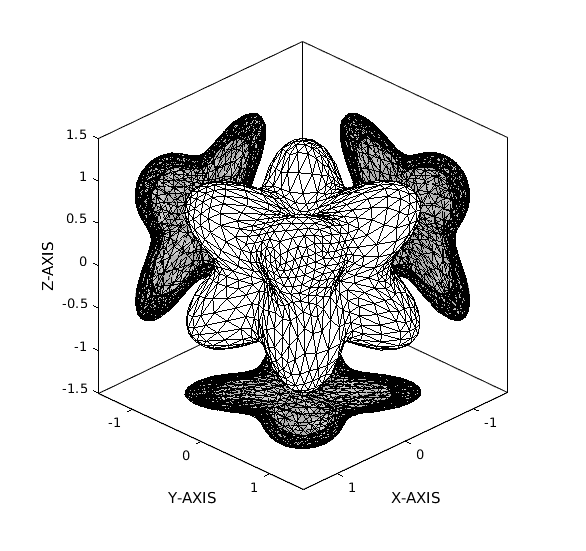}} 
	
	\subfloat[]{
		\label{fig:bumpy_cube_exact}
		\includegraphics[width=5cm,height=5cm]{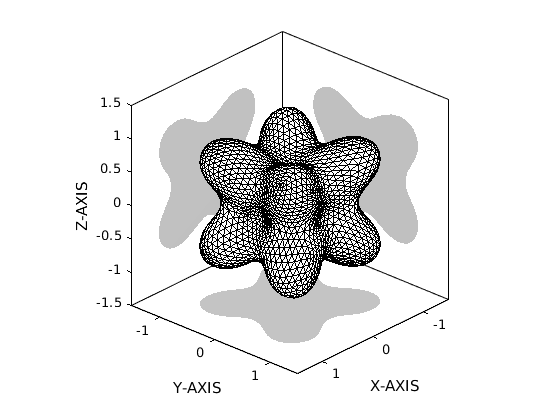}} 
	\caption{Reconstruction of a bumpy cube Fig.~\ref{fig:bumpy_cube_exact} from an initial guess (\ref{fig:bumpy_cube_initial}) after 19 iterations (\ref{fig:bumpy_cube_exact}) with $J(\beta) = 4.43 \cdot 10^{-2}$.}
	
\end{figure}

Next, we examine the performance of this algorithm on a more complex/challenging target illustrated in Fig.~\ref{fig:bumpy_cube_exact}. This geometry contains both convex and non-convex features which demonstrates the efficiency of MHs representation as well present an overall challenging reconstruction problem. In this case, we only use VSRM to obtain an initial guess. This is done using 16 incident planes waves in directions  $\mathbf{D}_{2} = \left \{ (\theta_2, \phi_2) \right \}$ and frequencies $\left \{ 100 Hz, 125 Hz, \cdots, 800 Hz\right \}$. The goal fields are generated for a set of 9 incident fields $\mathbf{D}_{1} = \left \{ (\theta_1, \phi_1) \right \}$ operating at $f_{1} = 400$ Hz. As such, the initial mesh is discretized with 1202 vertices and 2400 patches. The reconstruction routine was done using 40 MHs. The final result is shown in Fig.~\ref{fig:bumpy_cube_final}; the algorithm was  able to reconstruct the object and its features accurately from noisy data. Notably, the concave features of the object, i.e., the protruding lobes, were successfully recovered. Furthermore, the algorithm was able to recover the correct size and width of the scattering object. Again, we find their relative placement and generic dimensions agree very well with the target. 

\begin{table}[!h]
	\begin{center}
		\begin{tabular}{  m{3.0cm} | m{1.5cm} | m{1.5cm} | m{1.5cm}| m{1.5cm} |m{1.5cm}  } 
			\hline
			No. of iter. & Initial & 5 & 10 & 15 & 19 \\ 
			\hline
			$H\left(\Gamma_{G},\Gamma_{n} \right)$ & 2.48E-1 & 1.68E-1 & 9.96E-2 & 5.72E-2 & 4.18E-2 \\
			\hline
			$S_{err}\left(\Gamma_{G},\Gamma_{n}\right)$ & 1.15E-1 & 2.76E-3 & 1.90E-2 & 2.17E-2 & 1.41E-2 \\ 
			\hline
			$J(\beta_{n})$ & 3.21 & 1.63 & 1.45E-1 & 7.36E-2 & 4.44E-2 \\ 
			\hline
		\end{tabular}
	\end{center}
	\caption{Evolution of the error metrics with respect to iteration for the bumpy cube Ref~(\ref{fig:bumpy_cube_exact}) starting from Ref~(\ref{fig:bumpy_cube_initial}).}
	\label{table4}
\end{table}

\newpage

\subsection{Muli-resolution Reconstruction}

\begin{figure}[!h]
	\centering
	\subfloat[]{
		\label{fig:spot_initial}
		\includegraphics[width=5cm,height=5cm]{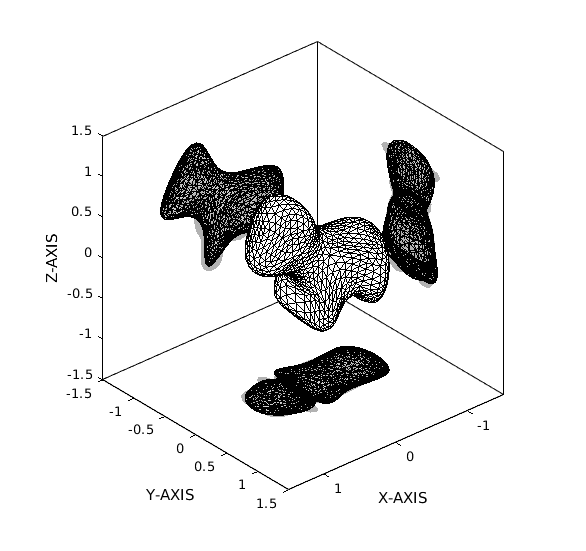} } 
	\hspace{1cm} 
	\subfloat[]{
		\label{fig:spot_final}
		\includegraphics[width=5cm,height=5cm]{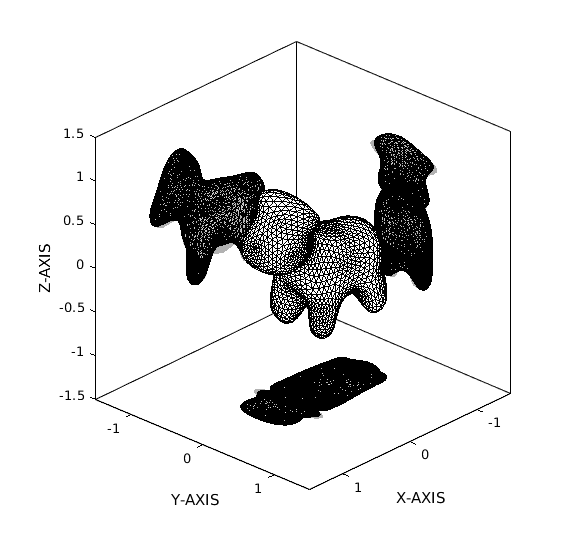} } 
	
	\subfloat[]{
		\label{fig:spot_exact}
		\includegraphics[width=5cm,height=5cm]{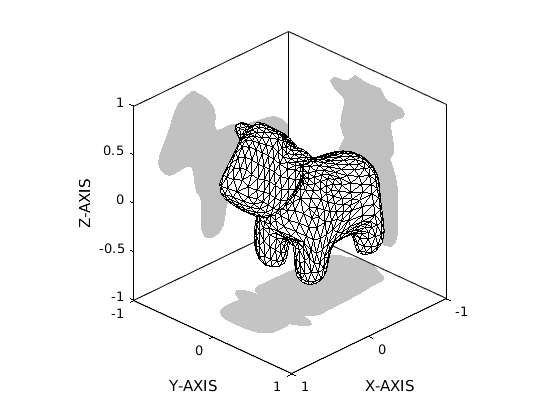} } 
	
	\caption{Reconstruction of Spot Fig.~\ref{fig:spot_exact} from an initial guess (\ref{fig:spot_initial}) after 63 iterations (\ref{fig:spot_exact}) with $J(\beta) = 2.44 \cdot 10^{-1}$. \label{fig:multiresolution}}
\end{figure}

In this final example, we consider the reconstruction of a cow that fits into a 1.7 m x 1.7 m x 1.0 m box, see Fig.~\ref{fig:spot_exact}. This geometry is multi-scale, implying there are fine-scale surface features as well as coarse features, which can be utilized for demonstrating the efficacy of the multi-resolution feature of our reconstruction scheme. Furthermore, given the multi-scale surface features, both the initial guess and the reconstruction routine require a richer set of data relative to the previous examples. In this case, the initial guess, see Fig.~\ref{fig:spot_initial}, is constructed from 16 incident plane waves in directions $\mathbf{D}_{2} = \left \{ (\theta_2, \phi_2) \right \}$, for the frequencies $\left \{ 100 Hz, 150 Hz, \cdots, 1200 Hz Hz\right \}$. The shape is optimized at two independent frequencies, 550 Hz and 900 Hz. Again, to maintain the prescribed number of samples over the mesh as $\lambda/7$, we discretize the mesh at 1102 vertices and 2200 faces for the first frequency and and then 2402 vertices and 4800 faces for the second frequency. Furthermore, at 550 Hz, we use 100 MHs and at 900 Hz we use 150 MHs. At each of these frequencies, we optimize in bands of 50 MHs. This implies that we choose $dim \left\{ \mathcal{B}_i \right \} = 50$ for any $i$.

\begin{table}[!h]
	\begin{center}
		\begin{tabular}{  m{3.0cm} | m{1.5cm} | m{1.5cm} | m{1.5cm}| m{1.5cm} |m{1.5cm}  } 
			\hline
			No. of iter. & Initial & 29 & 39 & 52 & 63 \\ 
			\hline
			$H\left(\Gamma_{G},\Gamma_{n} \right)$ & 1.20E-1 & 1.09E-1 & 1.01E-1 & 1.29E-1 & 9.00E-2 \\
			\hline
			$S_{err}\left(\Gamma_{G},\Gamma_{n}\right)$ & 3.41E-2 & 1.05E-2 & 5.84E-3 & 2.95E-2 & 7.87E-3 \\ 
			\hline
			$J(\beta_{n})$ & 6.69E-1 & 3.15E-2 & 1.93E-2 & 5.10E-1 & 2.44E-1 \\ 
			\hline
		\end{tabular}
	\end{center}
	\caption{Evolution of the error metrics with respect to iteration for the multiresolution cow Ref~(\ref{fig:spot_exact}) starting from Ref~(\ref{fig:spot_initial}).}
	\label{table5}
\end{table}

Reconstruction is done sequentially, in that we choose a threshold error for $\mathcal{B}_1$ and once that is reached, we optimize $\mathcal{B}_2$ and so on. Once we reach a local minima for the last band, we take our current optimized shape and refined it, such that we can optimize for the next frequency at the proper sampling rate and so on. As is evident from Fig.~\ref{fig:spot_final}, the algorithm recovers the major significant features of the object. While it does recover all features, locations and shape, it does not completely capture the legs or ear to high fidelity. This is largely due to reconstruction at two frequencies only that are not sufficient to resolve rather small features.

\section{Summary}
In the paper, we proposed a novel optimization scheme for shape reconstruction; the crux of our contributions rely on using manifold harmonics as a foundation for shape reconstruction. We have exploited the compression and multi-resolution framework that these harmonics provide to address every facet of the reconstruction process; from representing geometry to physics on the manifold. These harmonics rely on an underlying isogeometric analysis framework built on subdivision basis sets. A number of results demonstrate the viability and efficiency of the proposed approach on a set of challenging targets. As we proceed along this line of research, we anticipate developing a number of ideas that exploit manifold harmonics for inverse design; these include developing maps directly between geometry and far fields, multi-resolution editing, integration with machine learning, and so on. Several of these papers are underway, will be appear in the literature shortly. 
\section*{Acknowledgments}
The authors would acknowledge computing support from the HPC Center at Michigan State University and financial support from NSF via CMMI-1725278. Special thanks to Emeritus Professor Alex Diaz for all of the support and guidance he provided.




 \bibliographystyle{elsarticle-num-names} 
 \bibliography{references}


\end{document}